\documentclass[11pt]{article}
\usepackage{amsmath}
\usepackage{amssymb}
\usepackage{amsthm}
\usepackage{amsfonts}
\usepackage{graphicx}
\usepackage{pdfpages}
\usepackage{xcolor,colortbl}
\usepackage{multicol}
\usepackage{url}
\usepackage{subfigure}
\usepackage{authblk}  % Package for author and affiliation formatting

\usepackage{mathtools}

\DeclarePairedDelimiter\floor{\lfloor}{\rfloor}

\newcommand{\cF}{\mathcal{F}}

\def\P{{\mathbb P}}

\definecolor{Gray}{gray}{0.92}
\definecolor{LightCyan}{rgb}{0.88,1,1}

\newcolumntype{a}{>{\columncolor{Gray}}c}
\newcolumntype{b}{>{\columncolor{white}}c}

\usepackage{bbm,epsfig,graphics,epic,color,rotating,color}
    \usepackage[left=2.1cm,right=2.1cm,top=2.3cm,bottom=2.3cm]{geometry}

\numberwithin{equation}{section}
\def\E{{\mathbb E}}
\newcommand{\dE}{\mathbb{E}}

\def\P{{\mathbb P}}

\newtheorem*{theorem*}{Theorem}
\newtheorem{theorem}{Theorem}[section]
\newtheorem{lemma}[theorem]{Lemma}
\newtheorem{proposition}[theorem]{Proposition}

\theoremstyle{definition}

\title{A stochastic model for the diffusion of competing opinions with trend-following, opposition, and indifference.}
\date{}
\author{Manuel González-Navarrete\thanks{Departamento de Matem\'atica y Estad\'istica, Universidad de La Frontera, Temuco, Chile. email: manuel.gonzaleznavarrete@ufrontera.cl}}
\begin{document}

\maketitle %
%\vspace{1cm} %
\thispagestyle{empty} %
\baselineskip=14pt

\begin{abstract}
We study a stochastic model for the diffusion of competing opinions in a population composed of three types of agents: trend-followers, opposers, and indifferent individuals. The decision dynamics are driven by reinforcement mechanisms, modulated by a latent trend process, allowing us to capture realistic features such as amplification, resistance, and randomness in opinion formation. We derive explicit formulas for the finite-time moments of the opinion count vector and establish a set of strong asymptotic results, including laws of large numbers, central limit theorems, laws of the iterated logarithm, and almost sure convergence of empirical distributions. In particular, we show how early fluctuations can persist or vanish depending on the balance between reinforcement and opposition. Our analysis relies on martingale techniques and offers closed-form expressions for key quantities, providing both theoretical insights and tools for simulations or applications in social dynamics, marketing, or information diffusion.

%We study a stochastic model for the diffusion of opposing opinions ($\mathcal{A}$ and $\mathcal{B}$) in a population composed of three types of individuals: followers, opposers and indifferent agents. We derive explicit formulas for the moments of the random vector $(N_n,M_n)$, which represents the number of individuals who have adopted each opinion. In addition, we prove a collection of classical limit theorems describing the fluctuations of this vector around its deterministic limit. As a complement, we also provide a long-term characterization of the average of the difference process $(N_n-M_n)$.
% A key parameter $\theta$, aggregating the effects of social influence, determines the system's asymptotic regime: subcritical ($\theta < 1/2$), critical ($\theta = 1/2$), or supercritical ($\theta > 1/2$), each displaying qualitatively distinct behaviors.
\end{abstract}
%
%\bigskip

%%%%%%%%%%%%%%%%%%%%%%%%%%%%%%%%%%%%%%%%%%%%%%
%%%%%%%%%%%%%%%%%%%%%%%%%%%%%%%%%%%%%%%%%%%%%%
%%%%%%%%%%%%%%%%%%%%%%%%%%%%%%%%%%%%%%%%%%%%%%
%%%%%%%%%%%%%%%%%%%%%%%%%%%%%%%%%%%%%%%%%%%%%%

\section{Introduction}
\label{sec:intro}

The spread of competing opinions can be related to political opinion formation, market behavior (e.g., consumer choices), cultural dissemination or to analyze the influence of mass media and social pressure. In general, this phenomenon is studied through the use of mathematical models, which are naturally assumed to have stochastic dynamics. Within this context, principles from statistical mechanics can provide valuable insights, as presented in \cite{Bou,FPTT,Gomez,Hegselmann,MOS}. Models like the voter model, Sznajd model, or multilayer network frameworks have been widely used for this purpose \cite{NySz,Szna,Peng2021}. These systems are interesting given their nonlinear behaviors. A central feature commonly observed is the emergence of phase transitions, where small changes in parameters induces significant changes in global opinion. Including convergence to a dominating opinion, or phenomena like bistability o polarization \cite{Ba,BHW2}.

Following this line of research, Bendor et al. \cite{BHW} introduced a model for the diffusion of opinions as a probabilistic choice process. This model has its roots on Festinger's hypothesis \cite{Fe}. A mathematical characterization was proposed in \cite{GL2}. In that work, the authors introduced a generalized model, whose main feature is that each individual in the population has one of the following three intrinsic characteristics: \emph{trend-follower}, who tend to adopt the dominant or more visible opinion; \emph{opponents}, who actively resist prevailing trends and promote alternative views; or \emph{indifferent}, who do not engage significantly and may serve as noise in the diffusion process.

The model is defined as follows. Let denote the alternatives, $\mathcal{A}$ and $\mathcal{B}$. Assume that there is a positive initial number of ``converted" decision makers, denoted $N_0$ and $M_0$. In every period, one decision maker makes up his (her) mind about whether to adopt $\mathcal{A}$ or $\mathcal{B}$. The decision taken at the $n$-step is denoted by a Bernoulli variable $X_n$. In this sense, $X_n=1$ means $n$-decision was the innovation $\mathcal{A}$ (and $X_n=0$ means $\mathcal{B}$-decision). Let us denote by $N_n= N_0 + X_1 + \ldots + X_n$ the number of decision makers which adopted $\mathcal{A}$ and $M_n$ those that adopted $\mathcal{B}$, before the $(n+1)$-th decision.

The model in \cite{GL2} assumed the existence of a latent trend process $\{Y_n\}_{n \in \mathbb{N}}$, which represents a characteristic being independently imputed to each individual and defined as follows. For $\alpha, \beta \in (0,1)$ such that $0\leq \alpha+\beta\leq1$, the trend process $\{Y_n\}_{n \in \mathbb{N}}$ is an i.i.d. sequence of ternary random variables (which is independent of $N_n$ and $M_n$) such that, for each $n \ge 1$
$$
\P(Y_{n} = y) = \left\{\begin{array}{cccc}
                     \alpha   & \mbox{if} & y= +1 & \mbox{(Trend-follower),}\\ 
                     \beta   & \mbox{if} & y=-1 & \mbox{(Against-trend),}\\
                       1 - \alpha - \beta & \mbox{if} & y=0 & \mbox{(Indifferent).}
                       
\end{array}
\right.
$$

In this sense, the random-trend diffusion model is a stochastic process $\{X_n\}_{n \in \mathbb{N}}$ taking values on $\{0,1\}^{\mathbb{N}}$ defined by the following conditional probabilities
\begin{equation}
\label{Mine}
\P(X_{n+1} = 1 | \mathcal{F}_n, Y_n)=a + bY_n\frac{N_n}{N_n+M_n},
\end{equation}
where $N_n$ (respec. $M_n$) is given by the number of $\mathcal{A}$-decisions (respec. $\mathcal{B}$-decisions) taken until step $n$, $\mathcal{F}_n = \mathcal{F}_n(N_n,M_n)$ is the filtration related to random variables $N_n, M_n$. We assume that $Y_n$ is independent of $\mathcal{F}_n$, for all $n \ge 1$. Parameters $a$ and $b$ satisfy $0 \le a + b \le 1$, and also the following condition: if $\beta\neq 0$, then $b\leq a$.

We remark, as exposed in \cite{BHW}, that each decision is taken by a linear combination of own and social influences, being represented by parameter $a$ and proportion ${N_n}/{(N_n+M_n)}$, respectively. Furthermore, this dynamics can be incorporated to more complex processes with competing opinions (or decisions), such that in \cite{Kim,Peng2021}.

In \cite{GL2} the asymptotic characterization based on Janson \cite{Janson} were provided, including the following results, where the notation $\xrightarrow{d}$ denominates convergence in distribution, and $D([0,\infty[)$ being the Skorokhod space of right-continuous functions with left-hand limits.

\begin{theorem} \label{continuous} (Gonz\'alez-Navarrete \& Lambert, 2019 \cite{GL2})
Consider the random-trend diffusion model. If $0 \le a+b \le 1$ and $b\le a$ we get

\begin{equation}
\label{teo1}
\lim_{n \to \infty}\left(\frac{N_n}{N_n+M_n},\frac{M_n}{N_n+M_n}\right) = \left(\frac{a}{1-\theta},1-\frac{a}{1-\theta}\right)  \ \ \mbox{a.s.},
\end{equation}
where $\theta=b(\alpha-\beta)$. In addition,
\begin{itemize}
\item[(i)]  If $\theta <  \frac{1}{2}$ then, for $n \to \infty$, in $D[0, \infty)$

\begin{equation}
\label{teo2i}
\frac{1}{\sqrt{n}}\left[(N_{\floor*{tn}},M_{\floor*{tn}}) - tn\left(\frac{a}{1-\theta},1-\frac{a}{1-\theta}\right) \right] \xrightarrow{d} W_t,
\end{equation}
where $W_t$ is a continuous bivariate Gaussian process with $W_0=(0,0)$, $\E(W_t)=(0,0)$ and, for $0 < s \le t$,

\begin{equation*}
\E(W_s W_t^T) = s\left(\frac{t}{s}\right)^{\theta} \dfrac{a(1-a-\theta)}{(1-\theta)^2(1-2\theta)} \left( 
\begin{array}{cc}
1
& 
-1
\\
-1
& 
1
\end{array}
\right) \ .
\end{equation*}

\item[(ii)]  If $\theta =  \frac{1}{2}$ then, for $n \to \infty$, in $D[0, \infty)$

\begin{equation}
\label{teo2ii}
\frac{1}{\sqrt{n^t\log(n)}}\left[(N_{\floor*{n^t}},M_{\floor*{n^t}}) - n^t\left(2a, 1-2a\right) \right] \xrightarrow{d} W_t,
\end{equation}
where $W_t$ is a continuous bivariate Gaussian process with $W_0=(0,0)$, $\E(W_t)=(0,0)$ and, \break for $0 < s \le t$,

\begin{equation*}
\E(W_s W_t^T) = 2sa(1-2a)\left(\begin{array}{cc}
1 & -1 \\
-1 & 1 
\end{array}
\right) \ .
\end{equation*}
%with $\mathbf{1} = (1 \ 1)$.

\item[(iii)]\label{superdiffusive}
Suppose that $ \frac{1}{2}< \theta < 1$. Then we have tightness for the sequence of two-dimensional random variables $\{V_n\}_{n\geq2}$ given by $
V_n := \frac{1}{n^{b(\alpha - \beta)}}\left[{(N_n, M_n) - n\left(\frac{a}{1-\theta},1-\frac{a}{1-\theta}\right)}\right]$. Furthermore, there exists {a} random vector $W = ({W_1}, \ {W_2})$ such that, as $n$ diverges
\begin{equation}
\label{teo2iii}
V_n  \longrightarrow W \ \ \  \mbox{a.s.}
\end{equation}

\end{itemize}
\end{theorem}

To study the stochastic fluctuations of the proportion processes around their deterministic limits in \eqref{teo1}, we define
%In the sequel we use the following sequences of random variables, denoting the fluctuations of the proportions of opinions $\mathcal{A}$ and $\mathcal{B}$, with respect to limiting quantities, given by
\begin{equation}
    \label{Nhat}
    \hat{N}_n:=\frac{N_n}{T_n} - \frac{a}{1-\theta},
\end{equation}
and 
\begin{equation}
    \label{Mhat}
  \hat{M}_n:=\frac{M_n}{T_n} - \left(1-\frac{a}{1-\theta}\right),
\end{equation}
where $T_n = N_n + M_n$. The vector $(\hat{N}_n, \hat{M}_n)$ captures the deviations of the adoption proportions from their limit values.
We aim to derive the asymptotic characterization of this random vector. In addition, we will study the average cumulative difference process (as introduced in \cite{Grill}), which is defined by:
\begin{equation}
\label{centerofmass}
C_n = \frac{1}{n}\sum_{i=1}^n \left( N_i - M_i\right),
\end{equation}
which measures the average imbalance in opinion adoption over time.

In this work we complement Theorem \ref{continuous} by providing an explicit formula for the moments of random vector $(N_n,M_n)$ at finite time. Posteriorly, we state a collection of limit theorems for the behavior of the fluctuations in \eqref{Nhat} and \eqref{Mhat}, of particular interest is the characterization of the limit random vector $W$ in \eqref{teo2iii}. Moreover, we include a long-term characterization of the average difference process in \eqref{centerofmass}. The main tool is the novel martingale approach recently developed in \cite{BV2021,GLV2024,GLV2025,gue}. Our results show how the initial distribution of opinions influences the limit quantities.

The rest of the paper is organized as follows: in the next section, we provide the main results. The proofs of the Proposition and Theorems are presented in Section \ref{se:proofs}. Section \ref{se:remarks} contains some concluding remarks about the interpretation of our results and their mathematical significance. We conclude with some technical lemmas included in the appendix.

\section{Main Results}
\label{se:results}

This section presents the main results of this work. Let denote $T_n:=N_n+M_n = n +T_0$, and for $x \in \mathbb{R}$, the Pochhammer symbol $(x)_n:=\frac{\Gamma(n+x)}{\Gamma(x)}$, where $n$ is a non-negative integer, and $\Gamma$ is the gamma function. Then we obtain that

\begin{proposition}
\label{propo1}
    For any $n\ge 1$, we have
\begin{equation}
\label{Evector}
 \dE[(N_n,M_n)] = \dE(N_n)(1,-1)+T_n(0,1)
\end{equation}
and
\begin{equation}
\label{2doMvector}
 \dE[(N_n,M_n)^T(N_n,M_n)] = \dE(N_n^2)\left(\begin{array}{cc}
1 & -1 \\
-1 & 1 
\end{array}
\right)+T_n\dE(N_n)\left(\begin{array}{cc}
0 & 1 \\
1 & -2
\end{array}
\right)
+T_n^2\left(\begin{array}{cc}
0 & 0 \\
0 & 1
\end{array}
\right),
\end{equation}
where
\begin{equation}
\label{lENn}
 \dE[N_n] = \frac{(\theta+T_0)_n}{(T_0)_n}\left(N_0-\frac{a}{1-\theta} T_0\right)+\frac{a}{1-\theta}T_n,
\end{equation}
and
\begin{eqnarray}
\label{ENnsec}
 \dE[N_n^2]  = \frac{(2\theta+T_0)_n}{(T_0)_n}\Bigg[ N_0^2+  \frac{1}{(2\theta+T_0)}\left( N_0 - \frac{a}{1-\theta}T_0\right) \left[ (2a-1)T_0+(\theta+T_0)\left(\frac{2a(2\theta+T_0)}{(\theta-1)} \right) \right] \notag \\
+  \frac{a(1-2a)}{(1-\theta)(2\theta-1)} \left(T_0 -  \frac{(T_0)_{n+1}}{(2\theta+T_0)_n}\right)- \frac{a^2}{(1-\theta)^2} \frac{\Gamma(T_0+2)}{\Gamma(T_0)} \left(1-\frac{(T_0+2)_{n-1}}{(2\theta+T_0)_{n-1}} \right)  \notag\\
- \frac{(\theta+T_0)}{(2\theta+T_0)}\left( N_0 - \frac{a}{1-\theta}T_0\right) \left[ (2a-1)\frac{(\theta+T_0+1)_{n-1}}{(2\theta+T_0+1)_{n-1}}+ \frac{2a}{\theta-1}\frac{(\theta+T_0+1)_{n}}{(2\theta+T_0+1)_{n-1}}\right]  \Bigg].
\end{eqnarray}
\end{proposition}

These results can be complemented by recursive calculations using the formula
\begin{equation}\label{EN^k}
\mathbb{E}[N_{n+1}^k]=\mathbb{E}[N_n^k]\left(1+k\frac{\theta}{T_n}\right)+\sum_{j=1}^{k-1} \left(a{k\choose j}+\frac{\theta}{T_n}{k\choose j+1}\right)\mathbb{E}[N_n^{k-j}] + a,
\end{equation}
 for $k\ge 2$.
 
 Now, we define a quantity that will be related to the variance of limit random vectors. Namely
\begin{equation}
\label{sigmaa}
\sigma^2_{a,\theta} = \dfrac{a(1-a-\theta)}{(1-\theta)^2(1-2\theta)} .
\end{equation}
The next results are classified in the three different regions stated in Theorem \ref{continuous}. At the $\theta < 1/2$ setting, we provide the following limit theorem:
\begin{theorem}
\label{T_conv}
Assume $\theta < 1/2$. Then, as $n \to \infty$, the following results hold:

\begin{itemize}
\item[(i)] (Almost sure convergence of empirical measures)  
\begin{equation}
\label{asclt}
\displaystyle\lim_{n\to\infty}\displaystyle\frac{1}{\log n}\sum_{j=1}^n\frac{1}{j}\left(\mathbb{I}\left\{\sqrt{j}\hat{N}_j \le x \right\},\mathbb{I}\left\{\sqrt{j}\hat{M}_j \le y \right\}  \right) = G (x,y) \quad \text{almost surely,}
\end{equation}
where $\mathbb{I}\{\cdot\}$ denotes the indicator function of events and $G(x,y)$ is the distribution function of a bivariate normal vector $(X, Y)$ with mean and covariance matrix
\begin{equation}
\label{meanG1}
\mathbb{E}[(X,Y)] = (0,0), \qquad
\mathbb{V}[(X,Y)] = \sigma^2_{a,\theta} \begin{pmatrix}
1 & -1 \\
-1 & 1
\end{pmatrix},
\end{equation}
respectively.
\item[(ii)] (Almost sure convergence of moments)  
For all integers $m \geq 1$, we have
\begin{equation}
\label{moments1}
\frac{1}{\log n} \sum_{j=1}^n j^{m-1} \left( \hat{N}_j^{2m}, \hat{M}_j^{2m} \right) \longrightarrow \frac{\sigma_{a,\theta}^{2m}(2m)!}{2^m m!} (1,1) \quad \text{almost surely.}
\end{equation}

\item[(iii)] (Law of the iterated logarithm)  
\begin{equation}
\label{LIL1}
\limsup_{n \to \infty} \pm \frac{ \sqrt{n} (\hat{N}_n, \hat{M}_n) }{ \sqrt{2 \log \log n} } = \sigma_{a,\theta} (1,1) \quad \text{almost surely.}
\end{equation}

\item[(iv)] (Asymptotic behavior of the average difference process)  
We have the almost sure convergence
\[
\lim_{n \to \infty} \frac{C_n}{n} = \frac{a}{1 - \theta} - \frac{1}{2},
\]
and the following central limit theorem:
\begin{equation}
\label{CLT-CM}
\sqrt{n} \left( \frac{C_n}{n} - \left( \frac{a}{1 - \theta} - \frac{1}{2} \right) \right) \xrightarrow{d} \mathcal{N} \left( 0, \frac{8 \sigma^2_{a,\theta}}{3(2 - \theta)} \right).
\end{equation}

\end{itemize}
\end{theorem}

We should emphasize that, in this region, the influence of trend-followers and opposers is relatively weak. The population reaches a stable macroscopic distribution of opinions, and fluctuations are well-behaved. This regime corresponds to weak reinforcement or noise-dominated dynamics.

We continue by displaying the corresponding asymptotic analysis at the critical value $\theta=1/2$. 

\begin{theorem}
\label{criticalregion}

Assume $\theta=1/2$. Then, as $n \to \infty$, the following results hold:
\begin{itemize}

\item [(i)] (Almost sure convergence of empirical measures)
\begin{equation}
\label{asclt2}
\begin{array}{lll}
\displaystyle\lim_{n\to\infty}  \dfrac{1}{\log\log n}\displaystyle\sum_{j=2}^n\frac{1}{j\log j}\left(\mathbb{I}\left\{\sqrt{\frac{j}{\log j}}\hat{N}_j \le x \right\},\mathbb{I}\left\{\sqrt{\frac{j}{\log j}}\hat{M}_j \le y \right\}  \right) = G(x,y)    & \text{ a.s.},
\end{array}
\end{equation}
where $G$ is the distribution function of a normal random vector with
\begin{equation}\label{meanG2}
\mathbb{E}[(X,Y)]= (0,0)  \  \text{ and } \
\mathbb{V}[(X,Y)]= 2a(1-2a)  \left(\begin{array}{cc}
1 & -1 \\
-1 & 1 
\end{array}
\right).
\end{equation}

\item [(ii)] (Almost sure convergence of moments)  
For all integers $m \geq 1$, we have

\begin{equation}
\label{moments2}
\begin{array}{lll}
  \displaystyle\frac{1}{\log \log n}\sum_{j=2}^n \left(\frac{1}{j \log j}\right)^{m+1} j^{m-1} \left(\hat{N}_j^{2m},\hat{M}_j^{2m}\right) \rightarrow \frac{[a(1-2a)]^{m}(2m)!}{ m!}(1,1)   & \text{ a.s.},
\end{array}
\end{equation}

\item[(iii)] (Law of the iterated logarithm)
\begin{equation}
    \label{LIL2}
        \limsup_{n \to \infty} \pm \frac{\sqrt{n}(\hat{N}_n, \hat{M}_n) }{\sqrt{2 \log n \cdot \log \log ( \log n)}} = \sqrt{ 2a(1-2a)} (1,1) \ \text{ a.s. }
    \end{equation}

    \item [(iv)] (Asymptotic behavior of the average difference process)  
We have the almost sure convergence
\[\displaystyle\lim_{n \to \infty} \frac{C_n}{n}=\frac{4a-1}{2},\]
and the following central limit theorem:
\begin{equation}
\label{CLT-CM}
\sqrt{\frac{n}{\log n}}\left(\frac{C_n}{n}-\frac{4a-1}{2}\right) \overset{d}{\longrightarrow} \mathcal{N} \left(0,\frac{32a(1-2a)}{9}\right).
\end{equation} 
\end{itemize}

  \end{theorem}
At this region, the process exhibits marginal reinforcement, that is, the dependence is strong enough to prevent standard Gaussian scaling ($n$ is substituted by $\log n$), but not strong enough to induce random long-term limits, as will be expressed in the next case, whose deductions are technically more delicate.

\begin{theorem}
\label{T-ASP-SR}
Assume $\theta> 1/2$. Then, as $n \to \infty$, the following results hold:

\begin{enumerate}

\item [(i)] (Strong convergence to a random vector)
\begin{equation}
\label{FSLLN-SR}
\lim\limits_{n\to\infty}  n^{1-\theta}\left(\hat{N}_n, \hat{M}_n\right) =  (W,-W) \ \text{a.s.} \ \text{and in } L^2.
\end{equation}
where $W$ is a non-degenerated random variable such that
\begin{equation*}\label{meanW}
\mathbb{E}[W]=\frac{\Gamma(T_0)}{\Gamma(\theta+T_0)}\left(N_0-\frac{a}{1-\theta}T_0\right),
\end{equation*}
and
\begin{eqnarray}
\label{varW}
\dE[W^2]&=&\frac{\Gamma(T_0)}{\Gamma(2\theta+T_0)}\Bigg[ N_0^2+ \frac{a}{1-\theta}T_0\left(\frac{(1-2a)}{(2\theta-1)}- \frac{a}{(1-\theta)}(T_0+1) \right)\notag \\
&+ &  \frac{1}{(2\theta+T_0)}\left( N_0 - \frac{a}{1-\theta}T_0\right) \left[ (2a-1)T_0+(\theta+T_0)\left(\frac{2a(2\theta+T_0)}{(\theta-1)} \right) \right] \notag \Bigg].
\end{eqnarray}

\item [(ii)] (The Gaussian fluctuations)
\begin{equation}
    \label{T_cnv_w}
    \sqrt{n^{2\theta-1}}\left( n^{1-\theta}\left(\hat{N}_n, \hat{M}_n\right)-(W,-W) \right) 
\overset{d}{\longrightarrow} \mathcal{N}  \left(0,\Sigma\right),
\end{equation}
where
$$\Sigma =\dfrac{a(1-a-\theta)}{(1-\theta)^2(2\theta-1)} \left(\begin{array}{cc}
1 & -1 \\
-1 & 1 
\end{array}
\right)$$
\item [(iii)] (Law of the iterated logarithms)
    \begin{equation}
    \label{LIL3}
        \limsup_{n \to \infty} \pm \frac{\sqrt{n^{2\theta-1}}\left( n^{1-\theta}\left(\hat{N}_n, \hat{M}_n\right)-(W,-W) \right)}{\sqrt{2 \log \log n}} = \sqrt{\dfrac{a(1-a-\theta)}{(1-\theta)^2(2\theta-1)}} (1,1)\ \text{ a.s. }
    \end{equation}

\item [(iv)] (Asymptotics of the average difference process)
    \begin{equation}
    \lim_{n \to \infty} n^{1 - \theta} \left( \frac{C_n}{n} -\left(\frac{a}{(1-\theta)}-\frac12\right) \right) = \frac{2W}{1 + \theta} \quad \text{a.s.}
    \end{equation}
    \end{enumerate}
\end{theorem}

In this region, the system exhibits path-dependence and reinforcement learning behavior. The random limit $W$ reflects the memory of initial fluctuations, suggesting the system has a kind of ``frozen" asymmetry, even though the model itself is symmetric in design. Moreover, once again, we emphasize that higher moments of $W$ can be obtained by employing \eqref{EN^k} recursively.

\section{Proofs}\label{se:proofs}

First, note that from \eqref{Nhat} and \eqref{Mhat}, by a straightforward calculation, we obtain $\hat{N}_n + \hat{M}_n =0$.
In other words, all the calculations needed to analyze the vector of fluctuations $\left(\frac{N_n}{T_n} - \frac{a}{1-\theta}, \frac{M_n}{T_n} -\frac{1-a-\theta}{1-\theta} \right)$, are equivalent to study the vector $\left(\hat{N}_n,-\hat{N}_n\right)$.

Therefore, we mainly analyze the one-dimensional process $(N_n)$. Let start with some basic properties about conditional expectation in this model. Denote by $\left(\mathcal{F}_n\right)$ the increasing sequence of $\sigma$-algebras $\mathcal{F}_n=\sigma\left(X_1,\ldots,X_n\right)$. First of all, note that by the law of total probability we obtain
\begin{equation}
\label{condprobfinal}
 \P(X_{n+1} = 1 | \mathcal{F}_n) = a + \theta \frac{N_n}{T_n},
\end{equation}
which implies the Bernoulli random variables $X_n$, satisfy,
\begin{equation}\label{expmov}
\mathbb{E}\left[X_{n+1}|\mathcal{F}_n\right]=\theta\frac{N_n}{T_n}+a \hspace{.5cm}\text{a.s.}
\end{equation}
Hence, given that $N_{n+1}=N_n +X_{n+1}$, we have
\begin{equation}\label{expos}
\mathbb{E}\left[N_{n+1}|\mathcal{F}_n\right]=\gamma_n N_n+a \hspace{.5cm}\text{a.s.,}
\end{equation}%
where $\gamma_n=1+\frac{\theta}{T_n}.$
In general, by using that $\mathbb{E}\left[X^k_{n+1}|\mathcal{F}_n\right]=\mathbb{E}\left[X_{n+1}|\mathcal{F}_n\right]$, for all integer $k\ge 1$, then
\begin{equation}\label{k_condi}
\mathbb{E}[N_{n+1}^k|\mathcal{F}_{n}]=N_n^k+\left(a+\theta\frac{N_n}{T_n}\right)\sum_{j=1}^k {k\choose j}N_n^{k-j}, \ \text{ a.s.}
\end{equation}
Therefore, if $k\ge 2$, we obtain
\begin{equation}\label{k_condi2}
\mathbb{E}[N_{n+1}^k|\mathcal{F}_{n}]=N_n^k\left(1+k\frac{\theta}{T_n}\right)+\sum_{j=1}^{k-1} \left(a{k\choose j}+\frac{\theta}{T_n}{k\choose j+1}\right)N_n^{k-j} + a, \ \text{ a.s.}
\end{equation}
Note that \eqref{EN^k} is obtained by applying expectation on this equation. Moreover, due to the recursive nature of this formula, we need to consider sequences defined over the integers $n,k\geq 1$ by
\begin{equation}\label{an}
a_{n,k}=\prod_{j=0}^{n-1}\frac{T_j}{T_j+k\theta}=\frac{\Gamma(n+T_0)\Gamma(k\theta+T_0)}{\Gamma(n+k\theta+T_0)\Gamma(T_0)}= \frac{(T_0)_n}{(k\theta+T_0)_n}\sim \frac{\Gamma(k\theta+T_0) }{\Gamma(T_0)}n^{-k\theta},
\end{equation}%
where we use the Pochhammer symbol. Let us use $a_n:=a_{n,1}$ to simplify notation, and define the sequence $(A_n)$, which is given by $A_0=0$ and for $n\geq 1$ as
\begin{equation}
   \label{An}
A_n=\sum_{j=1}^n a_j.
\end{equation}
To analyze the behavior of $A_n$, from Lemma $B.1$ of \cite{ERWBercu}, we know that for $b\neq a+1$,
\begin{equation}
\label{B.1Lemma}
\sum_{j=1}^{n-1}\frac{\Gamma(j+a)}{\Gamma(j+b)} = \frac{\Gamma(a+1)}{(b-a-1)\Gamma(b)} \left( 1- \frac{(a+1)_{n-1}}{(b)_{n-1}} \right).
\end{equation}
Therefore, by \eqref{an} and \eqref{An} we obtain
\begin{equation}
\label{B.1Lem}
\frac{A_n}{ a_nT_n}=\frac{1}{\theta-1}\left(\frac{T_0}{ a_nT_n}-1 \right),
\end{equation}
which implies
\begin{equation}
\left|\frac{A_n}{ a_nT_n}-\frac{1}{1-\theta}\right| \sim  \frac{1}{\theta-1}\frac{T_0}{ a_nT_n}.
\label{desiAn}
\end{equation}

\subsection{Proof of Proposition \ref{propo1}}

Now, note that from \eqref{expos} we have for $n=1,2,\ldots,$ that
%\eqref{lENn2} directly lets us see that
\begin{equation}
\label{lENn2}
  \dE[N_n]=  \frac{1}{ a_n}\left(N_0+a\cdot  \displaystyle\sum_{j=1}^{n}a_j\right) = \frac{1}{ a_n}\left(N_0+a\cdot A_n\right).
\end{equation}
Since $A_n=\frac{1}{\theta-1}\left(T_0- a_nT_n\right)$ from \eqref{B.1Lem}, and using the definition of $a_n$ in \eqref{an}, we obtain \eqref{lENn}.
 
In the study of the second moment, we apply expectation in relation \eqref{k_condi2} to obtain that 
\begin{equation}
\mathbb{E}[N_{n+1}^2]=\delta_n\mathbb{E}[N_{n}^2]+\varphi_n \mathbb{E}[N_n]+a,
\end{equation}
where, $\delta_n:=1+\frac{2\theta}{T_n}$, and $\varphi_n:=2a+\frac{\theta}{T_n}$, for $n \geq 1$.
However, from \eqref{lENn2} we may see that
\begin{eqnarray*}
\varphi_n \mathbb{E}[N_n] = \left( \frac{2a T_n+\theta}{a_nT_n}\right) \left(N_0+ a A_n   \right),
\end{eqnarray*}
for $n \geq 1$. Hence, it may be found recursively for $n \geq 1$, that
\begin{eqnarray}
\label{MSN4}
\mathbb{E}[N_{n}^2]&=&\frac{1}{a_{n,2}}\left(N^2_0+\sum_{j=0}^{n-1}a_{j+1,2}\left(a+ \E(N_j) \varphi_j \right)\right) \notag\\
&=& \frac{1}{a_{n,2}}\left(N^2_0+\sum_{j=0}^{n-1}a_{j+1,2}\left(a+  \left( \frac{2a T_j+\theta}{a_jT_j}\right) \left(N_0+ a A_j   \right) \right)\right), \notag
\end{eqnarray}
where $a_{n,2}$ as defined in \eqref{an}, and $a_0=1$ and $A_0=0$. Again, given that $A_n=\frac{1}{\theta-1}\left(T_0- a_nT_n\right)$, we have
\begin{eqnarray}
\label{MSN5}
\hspace{-5cm}
\E[N_n^2]  & = & \frac{1}{a_{n,2}}\Bigg(N^2_0+ \frac{a}{1-\theta}\sum_{j=0}^{n-1}a_{j+1,2} +  \frac{2a^2}{1-\theta} \sum_{j=0}^{n-1} a_{j+1,2} T_j   \hfill
\notag\\
&&  + \left[ N_0 - \frac{a}{1-\theta}T_0\right] \left( 2a\sum_{j=0}^{n-1}  \frac{ a_{j+1,2}}{a_j} + \theta  \sum_{j=0}^{n-1}\frac{ a_{j+1,2}}{a_jT_j}\right)\Bigg).\hfill 
\end{eqnarray}
Therefore, we need to analyze the following sums,
\begin{equation}
\label{Sum1}
\sum_{j=0}^{n-1}a_{j+1,2}=\sum_{j=1}^{n}a_{j,2} = \frac{1}{2\theta-1}\left(T_0- a_{n,2}T_n \right),
\end{equation}
where $a_{n,2}T_n=\frac{(T_0)_{n+1}}{(2\theta+T_0)_n}$. Moreover,
\begin{equation}
\label{Sum2}
\sum_{j=0}^{n-1}a_{j+1,2} T_j= \frac{\Gamma(2\theta+T_0)}{\Gamma(T_0)}   \sum_{j=1}^{n} \frac{\Gamma(j+T_0+1)}{\Gamma(j+2\theta+T_0)}  -  \sum_{j=1}^{n}a_{j,2},
\end{equation}
\begin{equation}
\label{Sum3}
 \sum_{j=0}^{n-1}  \frac{ a_{j+1,2}}{a_j} = a_{1,2} + \frac{\Gamma(2\theta+T_0)}{\Gamma(\theta+T_0)}   \left( \sum_{j=1}^{n-1} \frac{\Gamma(j+\theta+T_0+1)}{\Gamma(j+2\theta+T_0+1)}  -  \theta \sum_{j=1}^{n-1} \frac{\Gamma(j+\theta+T_0)}{\Gamma(j+2\theta+T_0+1)}\right),
 \end{equation}
and
\begin{equation}
\label{Sum4}
 \sum_{j=0}^{n-1}  \frac{ a_{j+1,2}}{a_jT_j} = \frac{a_{1,2}}{T_0} + \frac{\Gamma(2\theta+T_0)}{\Gamma(\theta+T_0)}  \sum_{j=1}^{n-1} \frac{\Gamma(j+\theta+T_0)}{\Gamma(j+2\theta+T_0+1)},
\end{equation}
this last sum also appeared in \eqref{Sum3}. Finally, using \eqref{B.1Lemma} recursively we complete the proof.

%Then, we obtain that
%$$\sum_{j=0}^{n-1}\frac{ a_{j+1,2}  A_j}{a_j} = \frac{1}{\theta-1}\sum_{j=0}^{n-1}\frac{ a_{j+1,2}  }{a_j}\left(T_0- a_jT_j\right) =\frac{T_0}{\theta-1}\sum_{j=0}^{n-1}\frac{ a_{j+1,2}  }{a_j}-\frac{1}{\theta-1}\sum_{j=0}^{n-1}a_{j+1,2} T_j$$
%and
%$$\sum_{j=0}^{n-1}\frac{ a_{j+1,2}  A_j}{a_jT_j}= \frac{T_0}{\theta-1}\sum_{j=0}^{n-1}\frac{ a_{j+1,2}  }{a_jT_j} -\frac{1}{\theta-1}\sum_{j=0}^{n-1}a_{j+1,2} $$

% In addition, \eqref{MSN4} lets us to observe that

\subsection*{The associated martingale to characterize the asymptotic behaviors.}

Now, we construct a martingale associated to $N_n$, that it, the sequence $(Z_n)$, given by $Z_0=0$ and for $n\geq 1$  by
\begin{equation} \label{martingale}
Z_n=a_{n,1} N_n-a A_n,
\end{equation}
where the sequence $(A_n)$ as defined in \eqref{An}. Then, we observe from \eqref{expos} and \eqref{martingale} that almost surely
\begin{eqnarray*}
\mathbb{E}[Z_{n+1}|\mathcal{F}_n]=a_{n+1}(\gamma_n N_n+a)-a A_{n+1}=a_n N_n-a A_n=Z_n.
\end{eqnarray*}
Thus, $(Z_n)$ is a discrete time martingale with respect to the filtration $\left(\mathcal{F}_n\right)$.
%Now, note that
%\begin{equation}
%\label{Snequal}
%\dE(N_n)= \displaystyle\frac{\dE(X_1)}{a_n} + \frac{a}{a_n} \sum_{j=1}^{n-1} a_{j+1}.
%\end{equation}
From the definition of the proposed martingale given in \eqref{martingale}, we observe that
\begin{eqnarray}\label{delta}
Z_n=\sum_{j=1}^n \Delta Z_{j} =\sum_{j=1}^n (Z_{j}-Z_{j-1})=\sum_{j=1}^n a_j \xi_j,
\end{eqnarray}
where $\xi_1=X_1 - \mathbb{E}[X_1]$ and by \eqref{expos} we get that, for all $n\geq 1$, $\xi_{n+1}=N_{n+1}-\mathbb{E}[N_{n+1}|\mathcal{F}_{n}]=N_{n+1}-(a+\gamma_{n} N_n)$.
Furthermore, note that almost surely $\xi_{n+1}=X_{n+1}-\E[X_{n+1}|\mathcal{F}_n]$, for all $k \ge 1$.
From this and \eqref{expmov}, we find for $k\geq 2$ that,
%%%%%\begin{equation}\label{k_condi}
%%%%%\mathbb{E}[N_{n+1}^k|\mathcal{F}_{n}]=N_n^k+\left(a+\theta\frac{N_n}{n}\right)\sum_{j=1}^k {k\choose j}N_n^{k-j}, \ \text{ a.s., }
%%%%%\end{equation}
%which in turn leads us to
\begin{equation}
%\label{M3EPS}
\dE[\xi_{n+1}^k|\cF_n]= \sum_{j=0}^{k-2} {k\choose j} (-1)^j \left(a+\theta\frac{N_n}{T_n}\right)^{j+1}+(-1)^{k-1}(k-1)\left(a+\theta\frac{N_n}{T_n}\right)^{k}
\hspace{1cm}\text{a.s.,}
\notag
\end{equation}
that implies that
\begin{equation} \label{momentosk}
\sup_{n\geq 0} \mathbb{E}\left[\xi^k_{n+1}|\mathcal{F}_n\right] < \infty, \ \text{ a.s.,}
\end{equation}
because $\E[X_{n+1} | \mathcal{F}_n] \leq 1$, almost surely. In particular, 
\begin{equation} \label{segunda}
\mathbb{E}\left[\xi^2_{n+1}|\mathcal{F}_n\right]=\left(a+\theta\frac{N_n}{T_n}\right)\left(1-\left(a+\theta\frac{N_n}{T_n}\right)\right),
\ \text{a.s.}
\end{equation}
Hence, by analyzing the polynomial $p(x):=x-x^2$, for $0\le  x \le 1$, we have that
\begin{equation} \label{cotasegundomom}
\sup_{n\geq 0} \mathbb{E}\left[\xi^2_{n+1}|\mathcal{F}_n\right] \leq 1/4,  \ \text{a.s.}
\end{equation}
On the same line, 
\begin{equation}
\label{M4EPS}
\dE[\xi_{n+1}^4|\cF_n]=\left( a+\theta\frac{N_n}{T_n} \right)-4\left( a+\theta\frac{N_n}{T_n} \right)^2+6\left( a+\theta\frac{N_n}{T_n} \right)^3 
-3\left( a+\theta\frac{N_n}{T_n} \right)^4,
\ \text{a.s.}
\notag
\end{equation}
By maximizing the polynomial $p(x)=x-4x^2+6x^3-3x^4$, we conclude that
\begin{equation}\label{cuarta}
\sup_{n\geq 0} \mathbb{E}\left[\xi^4_{n+1}|\mathcal{F}_n\right] \leq 1/12, \ \text{a.s.}
\end{equation}
At this point, we empathize that the proofs will be based on the papers \cite{BV2021,GLV2024,heyde1977central}. Then, we need to define the predictable quadratic variation of {$(Z_n)$}, which satisfies for all $n \geq 1$, that
\begin{equation} \label{procrec}
\langle Z \rangle_n = \sum_{j=1}^n \mathbb{E}[\Delta Z_j^2| \mathcal{F}_{j-1}] = O(v_n),
 \end{equation}
where
$$v_n=\sum_{j=1}^n a_j^2.$$
Then via standard results on the asymptotic of the gamma function, we conclude that, as $n$ goes to infinity, it holds that

\begin{enumerate}
\item If $\theta<1/2$ then 
\begin{equation}
\label{vn1}
\lim_{n \to \infty} \frac{v_n}{n^{1-2\theta}}= \frac{\Gamma^2(\theta+T_0)}{(1-2\theta)\Gamma^2(T_0)}.
\end{equation}

\item If $\theta=1/2$ then
\begin{equation}
\label{vn2}
\lim_{n \to \infty} \frac{v_n}{\log n} = \frac{\Gamma^2(T_0+\frac12)}{\Gamma^2(T_0)}.
\end{equation}
%Which of course, leads us to conclude that $(v_n)$ goes to infinity as fast as $(\log n)$.

\item If $\theta>1/2$ then, from \eqref{an} it is possible to deduce that $(v_n)$ converges into a finite value, more precisely
\begin{equation} \label{vn3}
\lim_{n \to \infty}v_n = \sum_{j=0}^ \infty \left(\frac{\Gamma(j+T_0)\Gamma(\theta+T_0)}{\Gamma(j+\theta+T_0)\Gamma(T_0)} \right)^2 = \setlength\arraycolsep{1pt}
{}_3 F_2\left(\begin{matrix}T_0& &T_0& &1\\&\theta+T_0&
&\theta+T_0&\end{matrix};1\right),
\end{equation} 
where the above limit is the generalized hypergeometric function.
\end{enumerate}

The particular region described in the theorems are directly related to these limits.

\subsection{Proof of Theorem \ref{T_conv}}
First, note that \eqref{teo1} from Theorem \ref{continuous} together with \eqref{sigmaa} and \eqref{segunda} imply the following almost sure convergence 
\begin{equation}
\label{varianzaxi2}
\lim_{n\rightarrow \infty} \mathbb{E}\left[\xi^2_{n+1}|\mathcal{F}_n\right] = \sigma^2_{a,\theta}(1-2\theta).
\end{equation}
 In order to demonstrate $(i)$, we employ Lemma \ref{Lema1}, then by using \eqref{cuarta} we get
\begin{eqnarray*}
\sum_{j=1}^{\infty}\frac{1}{v_j}\E\left[|\Delta Z_j|^2\mathbb{I}_{\{|\Delta Z_j|\geq\varepsilon\sqrt{v_j}\}}|\mathcal{F}_{j-1}\right]  
\leq  
\frac{1}{\varepsilon^2} \sum_{j=1}^{\infty}\frac{1}{v_j^2}\E\left[|\Delta Z_j|^4|\mathcal{F}_{j-1}\right]  \\
\leq  
\sup_{j\geq 1}\E[\xi_j^4| \mathcal{F}_{j-1}]\frac{1}{\varepsilon^2}\sum_{j=1}^{\infty}\frac{a_j^4}{v_j^2}  
\leq  
\frac{1}{12\varepsilon^2}\sum_{j=1}^{\infty}\frac{a_j^4}{v_j^2} 
\sim 
\frac{1}{12\varepsilon^2}\sum_{j=1}^{\infty}\frac{(1-2\theta)^2}{j^2} < \infty,
\end{eqnarray*}
which implies condition \eqref{cond1}. The second condition of Lemma \ref{Lema1} is analogously proved by using $a=2$. Therefore we get that
$$
\displaystyle\lim_{n\to\infty}\frac{1}{\log v_n}\sum_{j=1}^n\left(\frac{v_j-v_{j-1}}{v_j}\right)\mathbb{I}\left\{\frac{Z_j}{\sqrt{v_{j-1}}}\le x \right\} = G (x)  \ \mbox{a.s.,}
$$
where $G$ is the distribution function of the normal random variable $N(0,\sigma^2_{a,\theta}(1-2\theta))$. Now, since the explosion coefficient is given by $
f_j = \frac{v_j-v_{j-1}}{v_j}=\frac{a_j^2}{v_j} \sim \frac{1-2\theta}{j}$,
{and by observing that} $\log v_n \sim (1-2\theta)\log n$ {we obtain that}
\begin{equation}
 \label{Zjvj}
\frac{Z_j}{\sqrt{v_{j-1}}} \sim \sqrt{\frac{1-2\theta}{j}}\left(N_j-a\frac{A_j}{a_j}\right) \sim \sqrt{(1-2\theta)j}\left(\frac{N_j}{T_j}-a\frac{A_j}{a_jT_j}\right)   ,
\end{equation}
then, by \eqref{desiAn} and the fact that $a_nT_n\to\infty$, as $n\to\infty$, we complete the proof of part $(i)$ in Theorem \ref{T_conv}.

Now, we focus our attention on the part $(ii)$. The proof is based on Lemma \ref{Lema2}. Notice that we have already seen that  $f_n \to 0$, as {$n\rightarrow \infty$}, then by \eqref{Zjvj} and \eqref{momentosk} we complete the proof.

Furthermore, to prove $(iii)$, we first remark that \eqref{varianzaxi2} jointly with 
$$\sum_{j=1}^{\infty}\frac{a_j^4}{v_n^2}=\frac{\pi^2(1-2\theta)^2}{6}, \ ´
$$
and the law of iterated logarithm for martingales (see Lemma C.2 in \cite{bblil}), we get that
$$
\limsup_{n \to \infty}\frac{Z_n}{\sqrt{2\log\log n}} = -\liminf_{n \to \infty}\frac{Z_n}{\sqrt{2\log\log n}} = \sqrt{\sigma^2_{a,\theta}(1-2\theta)} \hspace{1cm}\text{a.s.,} 
$$
which implies that
$$
\limsup_{n \to \infty}\left(\frac{n}{2\log\log n}\right)^{1/2}\left(\frac{N_n}{T_n}-a\frac{ A_n}{T_na_n}\right) = \sqrt{\sigma^2_{a,\theta}} \hspace{1cm}\text{a.s.} 
$$
%Therefore, \eqref{B.1Lem} conduces us to
%$$\limsup_{n \to \infty}\left(\frac{n}{2\log\log n}\right)^{1/2}\left(\frac{N_n}{n}-p+\frac{\theta p}{na_n}\right) = \sqrt{\sigma^2_{a,\theta}} \hspace{1cm}\text{a.s.} $$
Note that, since $
\displaystyle\lim_{n\to \infty}\left(\frac{n}{2\log\log n}\right)^{1/2}\frac{1}{n^{1-\theta}} =0$, we get that the last term vanishes as $n$ diverges, which completes the proof.

Finally, for (iv) we use the strategy proposed in \cite{BV2021} by employing the law of large numbers established in \eqref{teo1} that states that, almost surely
\[\frac{D_n}{n}:=\frac{N_n-M_n}{n}\rightarrow \frac{2a}{1-\theta} -1, \]
since $D_n=2N_n-T_n$. In addition, using the Toeplitz lemma (see, for instance  \cite{Duflo}) we obtain that
\[\frac{\displaystyle\sum_{i=1}^n i \frac{D_i}{i}}{\displaystyle\sum_{i=1}^n i}=\frac{2\displaystyle\sum_{i=1}^n i \frac{D_i}{i}}{n(n+1)}\rightarrow \frac{2a}{1-\theta} -1\hspace{1cm} \text{a.s.,}\]
which by definition of \eqref{centerofmass} implies
\[\frac{C_n}{n}=\frac{1}{n^2}\sum_{i=1}^n D_i\rightarrow \frac{a}{(1-\theta)}  -\frac12 \hspace{1cm} \text{a.s.}\]
For proving \eqref{CLT-CM}, first of all, note that
\[
C_n=\int_0^1 D_{{\lfloor nt \rfloor}}dt=\int_0^1( 2 N_{{\lfloor nt \rfloor}} - \lfloor nt \rfloor -T_0)dt,
\]
and that $\frac{\lfloor nt \rfloor}{n}\sim t,$ hence
\[\frac{C_n}{n}=\int_0^1 \frac{D_{{\lfloor nt \rfloor}}}{n}dt \sim \int_0^1 \left(2 \frac{N_{{\lfloor nt \rfloor}}}{n} - t\right)dt,\]
implies that
\[\sqrt{n}\left(\frac{C_n}{n}+\frac{1}{2}-\frac{a}{(1-\theta)}\right) \sim 2\int_0^1 \frac{1}{\sqrt{n}} \left(N_{{\lfloor nt \rfloor}}-\frac{a}{(1-\theta)}nt\right) dt, \]
which, by \eqref{teo2i} converges weakly to \(\displaystyle2\int_0^1 W_t dt,\) whose distribution is $N(0,\nabla^2)$, where
\[\nabla^2 = \dE\Big[\Big(2 \int_0^1 W_t dt\Big)^2\Big]=4 \int_0^1 \int_0^t 2\dE[W_s W_t] ds dt = 8\sigma^2_{a,\theta} \int_0^1 \int_0^t s\left(\frac{t}{s}\right)^{\theta}  ds dt=\frac{8\sigma^2_{a,\theta} }{3(2-\theta)}.\]

\subsection{Proof of Theorem \ref{criticalregion}}

We will proceed in a similar way as in the proof of Theorem \ref{T_conv}. First, due to \eqref{an} and \eqref{vn2} we obtain that 
\begin{equation}\label{ll2}
\frac{a_j^4}{v_j^2}\sim \left(\frac{1}{n \log n}\right)^2,
\end{equation}
which implies that 
\begin{equation} \label{ll3}
\sum_{j=1}^\infty \frac{a_j^4}{v_j^2}<\infty.
\end{equation}
Now, in order to demonstrate (i) we note  that the conditions of Lemma \ref{Lema1} follows from \eqref{ll3}. In addition, it may be found from the definition of $Z_n$, the fact that $\left|\frac{A_n}{ a_nT_n}-2\right|\sim \frac{2\Gamma( T_0+1)}{\sqrt{n}\Gamma(\theta + T_0)}$, \eqref{vn2} and \eqref{ll2} that $\frac{Z_n}{\sqrt{v_{n-1}}} \sim \sqrt{\frac{n}{\log n}}\left( \frac{N_n}{T_n}-2a \right),$ which leads to 
$$
\displaystyle\lim_{n\to\infty} \displaystyle\frac{1}{\log \log n}\sum_{j=2}^n \frac{1}{j \log j} \mathbb{I}\left\{\sqrt{\frac{j}{\log j}}\hat{N}_j \le x \right\} = G(x)   \ \text{ a.s.,}
$$ where $G$ stands for distribution function of a $N(0,2a(1-2a))$ random variable.

\medskip

For (ii), we notice that condition of Lemma \ref{Lema2} holds in the same manner than in the $\theta<1/2$ regime. From \eqref{varianzaxi2} and \eqref{ll2}, we obtain that $f_n$ converges to zero as $n \rightarrow \infty$. Hence, we may conclude \eqref{moments2} from the definition of $(Z_n)$.

\medskip

To prove (iii), we use similar arguments as in the proof of Theorem \ref{T_conv}, based on \eqref{ll3}.

Finally, for item (iv), the almost sure convergence immediately
follows from \eqref{teo1} together with Toeplitz lemma. The proof of the asymptotic normality is similar to that of Theorem 2.6 in \cite{BL}.

%
%
%from the last part of Theorem 1.3.24 of \cite{Duflo} we have that $\frac{Z_n^2}{\log n}=O(\log \log n)\hspace{.2cm}\text{a.s.}$ Moreover, from the definition of $(Z_n)$ we may see that $\left|a_n N_n-a A_n\right|=O(\sqrt{\log n \log \log n}) \ \text{a.s}$. Hence, from \eqref{an} we know that $na_n^2\rightarrow \frac{\pi}{4}$, then
%\begin{equation} \label{llncrit}
%\left|\frac{N_n}{n}-\frac{a A_n}{na_n}\right|=O\left(\sqrt{\frac{\log n \log \log n}{n}}\right)\hspace{.2cm}\text{a.s.}
%\end{equation}
%Additionally, we observe from \eqref{an} that $A_n\sim \sqrt{n\pi},$ which implies that $\frac{A_n}{na_n}\rightarrow 2.$ Furthermore, it is posible to obtain that $\left|\frac{A_n}{n a_n}-2\right|\sim \frac{2}{\sqrt{n\pi}}.$ Which, together with \eqref{llncrit} and the triangle inequality the proof holds. 

\subsection{ Proof of Theorem \ref{T-ASP-SR}.}

%From decomposition \eqref{martingale2}, \eqref{epsn} and \eqref{segunda}, we have; for $n\geq 0$, that
For proving (i), we recall that $\sup_{n \geq 0} \mathbb{E}[\xi^2_{n+1}|\mathcal{F}_n] \leq 1/4$ almost surely and {$(v_n)$} is a non increasing sequence. Then, it follows from \eqref{vn3}  that
$$\sup_{n\geq 1} \mathbb{E}[Z_{n}^2]\leq \frac{1}{4}  \cdot \setlength\arraycolsep{1pt}
{}_3 F_2\left(\begin{matrix}T_0& &T_0& &1\\&\theta+T_0&
&\theta+T_0&\end{matrix};1\right) <\infty.$$
That is to say, martingale $(Z_n)$ is bounded in $L^2$. Thus, it converges in $L^2$ and almost surely to the random variable $Z=\sum_{j=1}^\infty a_j \xi_j.$ Now, from \eqref{B.1Lem} we have that
$$
\lim_{n \to \infty} a_nT_n\left( \frac{N_n}{T_n}-\frac{a}{1-\theta}\right) = Z-\frac{a}{1-\theta}T_0 \ \ \mbox{a.s.,}$$
then
\begin{equation}
    \label{Wexist}
\lim_{n \to \infty}n^{1-\theta}\left( \frac{N_n}{T_n}-\frac{a}{1-\theta}\right) = W :=\left(Z-\frac{a}{1-\theta}T_0\right) \frac{\Gamma(T_0)}{\Gamma(\theta+T_0)} \ \mbox{a.s.} 
\end{equation}

Moreover, given that $(Z_n)$ converges to $Z$ in $L^2$ we have that \eqref{B.1Lem} and definition of the limit random variable $W$ guide us to
$$\lim_{n\rightarrow \infty} \mathbb{E}\left[\left(n^{1-\theta}\left(\frac{N_n}{T_n}-\frac{a}{1-\theta}\right) -W\right)^2 \right]=0.$$
We will find now the first two moments of the limiting random variable $W$ given in \eqref{FSLLN-SR}. Therefore, \eqref{lENn2} implies $\mathbb{E}[Z_n]=N_0$, which by \eqref{Wexist} implies that
$$\mathbb{E}[W]=\frac{\Gamma(T_0)}{\Gamma(\theta+T_0)}\left(\mathbb{E} (Z)-\frac{a}{1-\theta}T_0\right)=\frac{\Gamma(T_0)}{\Gamma(\theta+T_0)}\left(N_0-\frac{a}{1-\theta}T_0\right).$$
We will find the second moment of the limiting random variable $L$. For this, note that
\begin{eqnarray}
\mathbb{E}[Z_n^2]&=&a_n^2\mathbb{E}[N_n^2]-2a a_n A_n\mathbb{E}[N_n]+a ^2A_n^2 \notag\\[0.1cm]
&=&a_n^2\mathbb{E}[N_n^2]-2a  A_n N_0-a ^2A_n^2. \label{mnsquaree}
\end{eqnarray}

%
%
%Then
%\begin{eqnarray}
%\hspace{-1cm}
%\dE[N_n^2]&=&\frac{\Gamma(n+2\theta)}{\Gamma(n)}\left[
%\frac{\alpha}{\Gamma(2\theta+1)} + p \left( (1-2a) \sum_{j=1}^{n-1}\frac{\Gamma(k+1)}{\Gamma(k+1+2\theta)} + 2 a  \sum_{j=1}^{n-1}\frac{\Gamma(k+2)}{\Gamma(k+1+2\theta)}\right)  \right\notag\\
%&+& \left \frac{(\alpha-p)}{\Gamma(\theta+1)} \left( \theta(1-2a)  \sum_{j=1}^{n-1}\frac{\Gamma(k+\theta)}{\Gamma(k+1+2\theta)}  + 2a  \sum_{j=1}^{n-1}\frac{\Gamma(k+\theta+1)}{\Gamma(k+1+2\theta)} \right) \right]
%\notag
%\end{eqnarray}
%where $b=(p-\alpha)(1-\theta)$.
Hence, \eqref{ENnsec}, \eqref{an} and \eqref{mnsquaree} lead us to
\begin{eqnarray}
\hspace{-1cm}
& & \dE[Z_n^2]  \sim \frac{\Gamma^2(\theta+T_0)}{\Gamma(T_0)\Gamma(2\theta+T_0)}\Bigg[ N_0^2+  \frac{1}{(2\theta+T_0)}\left( N_0 - \frac{a}{1-\theta}T_0\right) \left[ (2a-1)T_0+(\theta+T_0)\left(\frac{2a(2\theta+T_0)}{(\theta-1)} \right) \right] \notag \\
&+&  \frac{a(1-2a)}{(1-\theta)(2\theta-1)}T_0- \frac{a^2}{(1-\theta)^2} \frac{\Gamma(T_0+2)}{\Gamma(T_0)} \left(1-\frac{\Gamma(2\theta+T_0)}{\Gamma(T_0+2)} n^{2(1-\theta)} \right)  \notag\\
&-& \frac{(\theta+T_0)}{(2\theta+T_0)}\left( N_0 - \frac{a}{1-\theta}T_0\right) \left[ \frac{2a}{\theta-1}\frac{\Gamma(2\theta+T_0+1)}{\Gamma(\theta+T_0+1)}n^{1-\theta}\right]  \Bigg] - \frac{2a}{\theta-1}N_0T_0 \left[ 1-\frac{\Gamma(\theta+T_0)}{\Gamma(T_0+1)}n^{1-\theta}\right]  \notag\\
&-& \frac{a^2}{(\theta-1)^2}T_0^2 \left[ 1-2\frac{\Gamma(\theta+T_0)}{\Gamma(T_0+1)}n^{1-\theta}+\frac{\Gamma^2(\theta+T_0)}{\Gamma^2(T_0+1)}n^{2(1-\theta)}\right].
\end{eqnarray}
Finally,
\begin{eqnarray}
\label{limZn}
%\hspace{-1cm}
 \lim_{n\rightarrow \infty} \dE[Z_n^2]  &=&\frac{\Gamma^2(\theta+T_0)}{\Gamma(T_0)\Gamma(2\theta+T_0)}\Bigg[ N_0^2+  \frac{1}{(2\theta+T_0)}\left( N_0 - \frac{a}{1-\theta}T_0\right) \left[ (2a-1)T_0+(\theta+T_0)\left(\frac{2a(2\theta+T_0)}{(\theta-1)} \right) \right] \notag \\
&+ & \frac{a(1-2a)}{(1-\theta)(2\theta-1)}T_0- \frac{a^2}{(1-\theta)^2} \frac{\Gamma(T_0+2)}{\Gamma(T_0)}    \Bigg] - \frac{2a}{\theta-1}N_0T_0 - \frac{a^2}{(\theta-1)^2}T_0^2,
\end{eqnarray}
which implies that
\begin{eqnarray}
\dE[W^2]&=&\frac{\Gamma^2(T_0)}{\Gamma^2(\theta+T_0)}\left(\dE[Z^2]-\frac{2a}{1-\theta}T_0\dE[Z]+ \left(\frac{a}{1-\theta}T_0\right)^2\right)\notag\\
&= & \frac{\Gamma(T_0)}{\Gamma(2\theta+T_0)}\Bigg[ N_0^2+  \frac{1}{(2\theta+T_0)}\left( N_0 - \frac{a}{1-\theta}T_0\right) \left[ (2a-1)T_0+(\theta+T_0)\left(\frac{2a(2\theta+T_0)}{(\theta-1)} \right) \right]\notag \\
&+ & \frac{a}{1-\theta}T_0\left(\frac{(1-2a)}{(2\theta-1)}- \frac{a}{(1-\theta)}(T_0+1) \right) \notag \Bigg].
\end{eqnarray}
In what follows, we will demonstrate items (ii) and (iii).  This proof will be based on Lemma \ref{Heydelem}. In this sense, note that \eqref{varianzaxi2} and the bounded convergence theorem imply that $\displaystyle\sum_{j=1}^{\infty}\dE[(\Delta Z_j)^2]\sim \sigma^2_{a,\theta}(1-2\theta)\frac{\Gamma^2(\theta+T_0)}{\Gamma^2(T_0)}\displaystyle\sum_{j=1}^{\infty}\frac{1}{j^{2\theta}}$. Then, since $\theta> 1/2$, we have that $\displaystyle\sum_{j=1}^{\infty}\dE[(\Delta Z_j)^2] < \infty$.
Now, given \eqref{varianzaxi2}, we obtain that
\begin{equation*}
\begin{array}{ll}
\displaystyle\sum_{j=n}^{\infty}\dE[(\Delta Z_j)^2\vert\mathcal{F}_{j-1}]&\sim\displaystyle\sum_{j=n}^{\infty}  \sigma^2_{a,\theta}(1-2\theta) a_j^2
    \sim \sigma^2_{a,\theta}(1-2\theta)\frac{\Gamma^2(\theta+T_0)}{\Gamma^2(T_0)}\displaystyle\sum_{j=n}^{\infty}\frac{1}{j^{2\theta}}\\[0.4cm]
    &\sim\displaystyle\frac{ -\sigma^2_{a,\theta}\Gamma^2(\theta+T_0)}{n^{2\theta-1}\Gamma^2(T_0)}
    \sim \dfrac{a(1-a-\theta)}{(1-\theta)^2(2\theta-1)} n a_n^2 \  \  \text{ a.s. }
\end{array}
\end{equation*}
By using the bounded convergence theorem, it follows that
\begin{equation}
\label{r2conv}
    r_n^2 := \sum_{j=n}^{\infty}\dE[(\Delta Z_j)^2]\sim \dfrac{a(1-a-\theta)}{(1-\theta)^2(2\theta-1)}n a_n^2\  \  \text{ a.s. }
\end{equation}
Then, conditions a) and a') of Lemma \ref{Heydelem} are satisfied. Hence, \eqref{an}, \eqref{delta} and the expectation of \eqref{cuarta} guide us to obtain an upper bound for $\dE\left[\vert\Delta Z_{j+1}\vert^4\right]$, then, \eqref{r2conv} and \eqref{an} let us conclude that
\begin{equation*}
\begin{array}{c}
         \displaystyle\frac{1}{r_n^2}\sum_{j=n}^{\infty}\dE\left[(\Delta Z_{j+1})^2 \mathbb{I}_{\{\vert\Delta Z_{j+1}\vert\geq\varepsilon r_n\}}\right]  \leq \displaystyle\frac{1}{\varepsilon^2r_n^4}\displaystyle\sum_{j=n}^{\infty}\dE\left[\vert\Delta Z_{j+1}\vert^4\right] \\[0,4cm]
 \leq \displaystyle\frac{1}{12\varepsilon^2r_n^4}\sum_{j=n}^{\infty} a_j^4
         \sim \frac{1}{r_n^4}\displaystyle\sum_{j=n}^{\infty} \frac{1}{j^{4\theta}}\sim \displaystyle  n^{4\theta-2}n^{1-4\theta},
\end{array}
 \end{equation*}
which implies condition b) of Lemma \ref{Heydelem}. Then, by noticing that $Z_n-Z = a_n \left(N_n - np -n^{\theta}L\right)$, and using \eqref{r2conv}, the normality holds for the first coordinate. However, the convergence \eqref{T_cnv_w} can be obtained by using the Cramér-Wold theorem (see, for instance Theorem 29.4 in \cite{Bill}).

Additionally, for $\varepsilon>0$, and from similar arguments than above, we have that
\begin{equation*}
         \frac{1}{r_j}\dE\left[\vert\Delta Z_{j+1}\vert\mathbb{I}_{\{\vert\Delta Z_{j+1}\vert\geq\varepsilon r_j\}}\right] \leq\frac{1}{r_j} \frac{1}{\varepsilon^3r_j^3}\dE\left[\vert\Delta Z_{j+1}\vert^4\right]
         \leq \frac{a_j^4}{12 \varepsilon^3 r_j^4}
				\sim \frac{1}{j^2},
 \end{equation*}
that is, condition c) of the same Lemma is satisfied.
In addition, given that $\sum_{j=1}^{\infty} \frac{1}{r_j^4}\dE[(\Delta Z_j)^4]<\infty$, we obtain condition d).
Finally, let us denote the martingale difference $d_j:=\frac{1}{r_j^2} \left( (\Delta Z_j)^2 - \dE[ (\Delta Z_j)^2 \vert\mathcal{F}_{j-1}]\right)$, and observe that
\begin{equation*}
\begin{array}{ll}
         \displaystyle\sum_{j=1}^{\infty}\dE\left[d_j^2 \vert\mathcal{F}_{j-1}\right]  = \displaystyle\sum_{j=1}^{\infty}\frac{1}{r_j^4}\left(\dE[ (\Delta Z_j)^4 \vert\mathcal{F}_{j-1}]-\dE^2[ (\Delta Z_j)^2 \vert\mathcal{F}_{j-1}]\right)\\[0,4cm]
        \leq \displaystyle\sum_{j=1}^{\infty}\frac{1}{r_j^4}\dE[ (\Delta Z_j)^4 \vert\mathcal{F}_{j-1}]
         \leq \frac{1}{12}\displaystyle\sum_{j=1}^{\infty} \frac{a_j^4}{r_j^4}\leq \frac{1}{12}\left(\dfrac{(1-\theta)^2(2\theta-1)}{a(1-a-\theta)}\right)^2\displaystyle\sum_{j=1}^{\infty} \frac{1}{j^2} < + \infty.
\end{array}
 \end{equation*}
 Then, as a consequence of Theorem 2.15 in \cite{hall2014martingale}, we can use Corollary 2 from \cite{heyde1977central}, and therefore \eqref{LIL3} holds.

For item (iv), from the definition of \( C_n \) we have
\begin{equation}
n C_n = \sum_{j=1}^n D_j = \sum_{j=1}^n 2\left( \frac{N_j}{T_j} -\frac{a}{(1-\theta)}\right) T_j + \left(\frac{2a}{(1-\theta)}-1\right)\sum_{j=1}^n T_j.
\end{equation}

%which can be rewritten as
%\begin{equation} \label{eq:nCn}
%n C_n = \sum_{j=1}^n j^{1 - \theta} \left( \frac{D_j}{j} - \left(\frac{2a}{(1-\theta)}-1\right) \right) j^{\theta} +\left(\frac{2a}{(1-\theta)}-1\right)\cdot \frac{n(n+1)}{2}.
%\end{equation}

It follows from equation \eqref{FSLLN-SR} together with Toeplitz's lemma \cite{Duflo} that
\begin{equation} \label{eq:ToepCn}
\lim_{n \to \infty} \frac{1}{n^{1 + \theta}} \sum_{j=1}^n 2T_j\left( \frac{N_j}{T_j} -\frac{a}{(1-\theta)}\right) = \frac{2W}{1 + \theta} \quad \text{a.s.}
\end{equation}

Hence, we deduce that
\[
\lim_{n \to \infty} \frac{n}{n^{1 + \theta}} \left( C_n - \left(\frac{2a}{(1-\theta)}-1\right)\left(\frac{(n + 1)}{2} +T_0\right)\right) = \frac{2W}{1 + \theta} \quad \text{a.s.},
\]

which clearly leads to
\[
\lim_{n \to \infty} n^{1 - \theta} \left( \frac{C_n}{n} - \left(\frac{a}{(1-\theta)}-\frac12\right) \right) = \frac{2W}{1 + \theta} \quad \text{a.s.}
\]

\section{Concluding Remarks}
\label{se:remarks}

In this work, we analyzed a stochastic process for the diffusion of competing opinions, incorporating three types of agents: trend-followers ($\alpha$), who amplify the dominant opinion and reinforce trends; opposers ($\beta$), who counteract prevailing tendencies and promote balance; and indifferent individuals ($1-\alpha-\beta$), who introduce randomness that may stabilize or destabilize the system depending on their relative proportion.

Our finite-time characterization, particularly in Proposition~\ref{propo1}, established an explicit expression for the first and second moments of the opinion counts. This revealed a clear negative correlation between competing opinions and quantified the influence of initial conditions. These results are valuable for numerical simulations, model calibration, and predictions at moderate time scales, highlighting how parameters such as $(a,\theta,N_0,M_0)$ affect the system's evolution.

The key parameter $\theta = b(\alpha - \beta)$ condenses the interplay between reinforcement and opposition into a single quantity that governs the system’s long-term behavior. When opposition balances or outweighs reinforcement ($\theta < 1/2$), the system stabilizes and converges deterministically. At the critical point ($\theta = 1/2$), the system exhibits marginal reinforcement, where slight fluctuations or biases can lead to long-term asymmetries. In this case, the process becomes statistically fragile, and detecting stable patterns requires longer time series. When trend-following dominates ($\theta > 1/2$), early fluctuations persist and strongly influence the limiting distribution, breaking symmetry and yielding path-dependent outcomes.

These findings have broad implications for modeling social contagion, consumer behavior, political polarization, and information diffusion, where early randomness may have permanent consequences on collective dynamics.

From a mathematical perspective, we remarks that the use of martingale methods and explicit moment calculations offers sharp characterizations, going beyond mere existence results. The laws of the iterated logarithm and almost sure convergence of empirical measures, rarely available in models of this nature, provide strong asymptotic tools. Moreover, the explicit expression for the limit variable $W$ and its moments contributes to analytical tractability and facilitates simulations and applications.

An interesting direction for future work is to extend this framework to structured populations, such as multilayer or dynamic networks, where the influence of local neighborhoods, community structure, or temporal interactions may amplify or dampen global trends. In the mathematical context, it would also be valuable to explore Cramer moderate deviations \cite{FS} or the characterization of equilibrium times (i.e. $N_n=M_n$) \cite{Bertoin}. Even a multiple opinion model can be studied by employing similar techniques to those in \cite{BL,G2020,Qin,gue}.

%These extensions could further enhance the model’s applicability to real-world phenomena and provide deeper insights into the dynamics of opinion formation and diffusion.

%Another natural extension involves external influence or control mechanisms, such as targeted interventions or mass media effects, which can be incorporated through time-dependent parameters or exogenous shocks. 

\section*{Appendix}

In this section we include three Lemmas that are used in the proof of asymptotic results. In particular, the first one is an adaptation of Theorem A.1 from \cite{gue}, which was originally proved in \cite{cvgm}. The Lemma \ref{Lema2} is based on Theorem A.2 from \cite{gue}, its original version may be found in \cite{Cha}. Finally, Lemma \ref{Heydelem} may be obtained from Theorem 1 and Corollaries 1 and 2 from \cite{heyde1977central}.

Let a martingale $(Z_n)$ with $\Delta Z_n = Z_n - Z_{n-1} = a_n\xi_n$, where $(a_n)$ is  a deterministic real sequence and $\xi_n$ are random variables such that
$$\lim_{n\rightarrow \infty} \mathbb{E}\left[\xi^2_{n+1}|\mathcal{F}_n\right] = \sigma^2,$$
almost surely. Also denote $v_n=\sum_{j=1}^n a_j$, increasing to infinity.
Therefore,
\begin{lemma}
\label{Lema1}
Assume for all $\varepsilon > 0$ that
\begin{equation}
\label{cond1}
    \displaystyle\sum_{n=1}^{\infty} \frac{1}{v_n} \dE\left[\vert\Delta Z_n\vert^2 \mathbb{I}_{\{\vert\Delta Z_n\vert\geq\varepsilon\sqrt{v_n}\}}\vert\mathcal{F}_{n-1}\right] < \infty \ \  \text{ a.s.},
\end{equation}
and for some $a>0$ that,
\begin{equation}
\label{cond2}
    \displaystyle\sum_{n=1}^{\infty} \frac{1}{v_n^a} \dE\left[\vert\Delta Z_n\vert^{2a} \mathbb{I}_{\{\vert\Delta Z_n\vert\leq\sqrt{v_n}\}}\vert\mathcal{F}_{n-1}\right] < \infty \ \  \text{ a.s. }
\end{equation}
Then, as $n$ diverges, $(Z_n)$ satisfies that
\begin{equation}
\displaystyle\lim_{n\to\infty}\displaystyle\frac{1}{\log v_{n}}\sum_{k=1}^n \left(\frac{v_{k}-v_{k-1}}{v_{k}}\right) \mathbb{I}\left\{\frac{Z_k}{\sqrt{v_{k-1}} } < x \right\} = G (x) \ \  \text{ a.s.,}
\end{equation}
where $\mathbb{I}\{\cdot\}$ denotes  the indicator function of events and $G$ stands for the distribution function of a normal random variable $N(0,\sigma^2)$.
\end{lemma}

\begin{lemma}
\label{Lema2}
Assume that $(\xi_n)$ satisfies, for some integer $m\ge 1$ and for some real number $h > 2m$, that
\begin{equation}
\label{cond3}
    \displaystyle\sup_{n\ge 0}\dE\left[\vert\xi_{n+1}\vert^h\vert\mathcal{F}_{n}\right] < \infty \ \  \text{ a.s. }
\end{equation}
If in addition the explosion coefficient $f_n=\frac{v_{n}-v_{n-1}}{v_{n}}$ tends almost surely to zero, then
\begin{equation}
    \displaystyle\lim_{n\to\infty} \frac{1}{\log v_n}\sum_{k=1}^n f_k \left( \frac{Z_k^2}{v_{k-1}} \right)^m = \frac{\sigma^{2m} (2m)!}{2^m m!} \ \  \text{ a.s. }
\end{equation}
\end{lemma}

%
%%%%%%%%%%%%%%%%%%%%%%%%%%%%%%%%%%%%%%%%%%%%%%%%%%%%%%

%Finally, we state the following lemma without proof. However it 

\begin{lemma}
\label{Heydelem}
Suppose that $(Z_n)$ is a zero mean square-integrable martingale. Let $r_n^2:= \displaystyle\sum_{k=n}^{\infty} \dE[(\Delta Z_k)^2]$. If $\displaystyle\sum_{k=1}^{\infty} \dE[(\Delta Z_k)^2] < +\infty,$
and $\displaystyle\sum_{k=1}^{\infty}\frac{1}{r^2_k} \left[ (\Delta Z_k)^2 - \dE[ (\Delta Z_k)^2 \vert\mathcal{F}_{k-1}]\right] < +\infty ~ \text{ a.s.},$ we have the following:
\begin{itemize}
    \item [(i)] The limit $Z:= \sum_{k=1}^{\infty} \Delta Z_k $ exists almost surely and $Z_n \overset{L^2}{\rightarrow} Z$.
    \item [(ii)] Assume additionally that
%    \begin{fleqn}
		\begin{itemize}
        \item [a)] $\displaystyle\frac{1}{r_n^2}\sum_{k=n}^{\infty}\dE\left[(\Delta Z_k)^2\vert\mathcal{F}_{k-1}\right]\rightarrow 1$ as $n \to \infty$ in probability, and
        \item [b)] $\displaystyle\frac{1}{r_n^2}\displaystyle\sum_{k=n}^{\infty}\dE\left[(\Delta Z_{k+1})^2 \mathbb{I}_{\{\vert\Delta Z_{k+1}\vert\geq\varepsilon r_k\}}\right] \rightarrow 0$ for any $\varepsilon > 0$.
        
        Then, we have that
        \begin{equation*}
            \displaystyle\frac{Z-Z_n}{r_{n+1}}=\frac{\sum_{k=n+1}^{\infty}\Delta Z_k}{r_{n+1}} \overset{d}{\rightarrow}N(0,1).
        \end{equation*}
    \end{itemize}
%		\end{fleqn}
     \item [(iii)] If, in addition, the following conditions hold 
       \begin{itemize}
        \item [a')] $\displaystyle\frac{1}{r_n^2}\sum_{k=n}^{\infty}\dE\left[(\Delta Z_k)^2\vert\mathcal{F}_{k-1}\right]\rightarrow 1$ as $n \to \infty$ a.s.
        \item [c)] $\displaystyle\sum_{k=1}^{\infty}\frac{1}{r_k}\dE[ |\Delta Z_k| \mathbb{I}_{\{| \Delta Z_k|  > \varepsilon r_k\}}] < +\infty$ for any $\varepsilon > 0$, and
        \item [d)] $\displaystyle\sum_{k=1}^{\infty}\frac{1}{r_k^4}\dE[ (\Delta Z_k)^4 \mathbb{I}_{\{|\Delta Z_k|  \le \delta r_k\}}] < +\infty$ for some $\delta > 0$. 
    \end{itemize}
    Then, $\displaystyle\limsup_{n\to\infty}\pm \frac{Z-Z_n}{\sqrt{2\cdot r^2_{n+1} \log |\log r^2_{n+1} |}} = 1$ a.s.
\end{itemize}
\end{lemma}


\begin{thebibliography}{}

\bibitem{Ba}
{\sc Banerjee, A.V.} (1992)
A simple model of herd behavior.
{\em Quarterly Journal of Economics } \textbf{107}, 797-818.


\bibitem{BHW}
{\sc Bendor, J., Huberman, B. and Wu, F.} (2009)
Management fads, pedagogies, and other soft technologies.
{\em J. Econ. Behav. Organ.} \textbf{72}, 290-304.

\bibitem{bblil} {\sc Bercu, B.} (1998)
{\it Central Limit Theorem and Law of Iterated Logarithm for Least Squares Algorithms in Adaptive Tracking}. SIAM Journal on Control and Optimization.  \textbf{36}(3), 910--928.

\bibitem{cvgm} {\sc Bercu, B.} (2004) {\it On the convergence of moments in the almost sure central limit theorem for martingales with statistical applications.} Stochastic Process. Appl. \textbf{111}(1), 157-173.

\bibitem{ERWBercu}  {\sc Bercu, B.} (2018) {\it A martingale approach for the elephant random walk.} J. Phys. A: Math. Theor. \textbf{51}, 015201.

\bibitem{BL} {\sc Bercu, B. and Laulin, L.} (2021)
{\it On the center of mass of the elephant random walk.} Stoch. Process. Appl. \textbf{133}, 111--128

\bibitem{BV2021}
{\sc Bercu, B. and. Vázquez Guevara, V. H} (2022)
{\it Further results on the minimal random walk,} J. Phys. A: Math. Theor. \textbf{55}(41), 415001.


\bibitem{Bertoin} {\sc Bertoin, J.} (2022)
{\it Counting the zeros of an elephant random walk}. Trans. Amer. Math. Soc. 375, 5539-5560.


\bibitem{BHW2}
{\sc Bikhchandani, S., Hirshleifer, D. and Welch, I.} (1992)
A theory of fads, fashion, custom, and cultural change as informational cascades.
{\em Journal of Political Economy } \textbf{100}, 992-1026.

\bibitem{Bill}
{\sc Billingsley, P.} (2008)
Probability and measure. John Wiley \& Sons.


\bibitem{Bou}
{\sc Bouchaud, J.P.} (2013)
Crises and Collective Socio-Economic Phenomena:
Simple Models and Challenges.
{\em J. Stat. Phys.} \textbf{151}(3-4), 567--606.


\bibitem{Cha}{\sc Cha\^abane, F.}
\newblock {\it Version forte du theoreme de la limite centrale fonctionnel pour les martingales.}
C. R. Acad. Sci. Paris \textbf{323}: 195-198 (1996).



%\bibitem{Dr}
%{\sc Drezner, Z.} (2006)
%On the Limit of the Generalized Binomial Distribution.
%{\em Comm. Statist. Theory Methods} \textbf{35}(2), 209--221.

%\bibitem{DF}
%{\sc Drezner, Z. and Farnum, N.} (1993)
%A generalized binomial distribution.
%{\em Comm. Statist. Theory Methods} \textbf{22}, 3051--3063.

\bibitem{Duflo}{\sc Duflo, M.} (1997) {\it Random Iterative Models.} Springer Verlag, Berlin.
%
%\bibitem{DR} {\sc Durrett, R. and Resnick, S.} (1978)
%{\it Functional limit theorems for dependent random variables.} Ann. Probab. \textbf{6}: 829-846.

\bibitem{FS} {\sc Fan, X and Shao, Q-M}. (2024)
Cram\'er’s moderate deviations for martingales with applications. Ann. inst. Henri Poincare (B) Probab. Stat. \textbf{60}(3): 2046--2074.

\bibitem{Fe}
{\sc Festinger, L.} (1954)
A theory of social comparison processes.
{\em Hum. Relat.} \textbf{7}(2), 117--140. 

\bibitem{FPTT}
{\sc Furioli, G., Pulvirenti, A., Terraneo E. and Toscani, G.} (2017)
Fokker-Planck equations in the modelling of socio-economic phenomena. Math Mod Meth Appl Sci \textbf{27}(1): 115--158.

\bibitem{Gomez} {\sc Gómez-Serrano, J., Graham, C. and  Le Boudec,J.-Y.} (2012)
{\it The bounded confidence model of opinion dynamics,} Math. Models Methods Appl. Sci., \textbf{22}, 1--46.

\bibitem{G2020} {\sc Gonz\'alez-Navarrete, M.} (2020)
{\it Multidimensional walks with random tendency.}
J. Stat. Phys. \textbf{181}(3): 1138--1148.

%\bibitem{GL}
%{\sc Gonz\'alez-Navarrete, M. and Lambert, R.} (2018) Non-Markovian random walks with memory lapses. {\em J. Math. Phys.} \textbf{59}, 113301.

\bibitem{GL2} {\sc Gonz\'alez-Navarrete, M. and Lambert, R.} (2019)
{\it The diffusion of opposite opinions in a randomly biased environment.}
J. Math. Phys. \textbf{60}, 113301.


\bibitem{GLV2024} {\sc Gonz\'alez-Navarrete, M., Lambert, R. and V\'azquez-Guevara, V.H.} (2024)
{\it A complete characterization of a correlated Bernoulli process.}
Electron. Commun. Probab. \textbf{29}(69), 1-12.

\bibitem{GLV2025} {\sc Gonz\'alez-Navarrete, M., Lambert, R. and V\'azquez-Guevara, V.H.} (2025)
{\it On the asymptotic analysis of lazy reinforced random walks: a martingale approach.}
{J. Math. Anal. Appl.} \textbf{549}(2), 129520.

\bibitem{Grill}
 {\sc Grill, K.} (1988) On the average of a random walk, Statist. Probab. Lett. \textbf{6}(5), 357–361.

\bibitem{hall2014martingale} {\sc Hall, P. and Heyde, C.C.} (1980)
{\it Martingale Limit Theory and Its Application.} Academic Press, New York.


\bibitem{Hegselmann} {\sc Hegselmann, R. and Krause. U.} (2002)
Opinion dynamics and bounded confidence models, analysis, and simulation. J. Artif. Soc. Soc. Simul., 5(3).



%
%\bibitem[Hall and Heyde(1980)]{H&H}
%{\sc Hall, P. and Heyde, C.C.} (1980)
%{\em Martingale limit theory and its application.}
%Academic Press, New York.

\bibitem{heyde1977central} {\sc Heyde, C.C.} (1977)
\newblock {\it On central limit and iterated logarithm supplements to the martingale convergence theorem.}
J. Appl. Probab. \textbf{14}(4), 758-775.

%\bibitem{He}
%{\sc Heyde, C.C.} (2004)
%Asymptotics and criticality for a correlated Bernoulli process.
%{\em Aust. N. Z. J. Statist.} \textbf{46}, 53--57.

\bibitem{Janson}
{\sc Janson, S.} (2004)
Functional limit theorems for multitype branching processes and generalized P\'olya urns.
{\em Stoch. Process. Appl.} \textbf{110}, 177--245.


\bibitem{Kim}
{\sc Kim, S.-J., Naruse, M., Aono, M., Ohtsu, M. and Hara, M.} (2013)
Decision Maker Based on Nanoscale Photo-Excitation Transfer. Sci Rep. \textbf{3},
2370.

%\bibitem{LH}
%{\sc Liu, X. and Hong, Y.} (2015)
%Analysis of railroad tank car releases using a generalized binomial model.
%{\em  Accid. Anal. Prev.} \textbf{84}, 20--26.

%\bibitem{LH2}
%{\sc Liu, X., Liu, C. and Hong, Y.} (2017)
%Analysis of multiple tank car releases in train accidents.
%{\em  Accid. Anal. Prev.} \textbf{107}, 164--172.


\bibitem{MOS}
{\sc Moharir, S., Omanwar, A.S. and Sahasrabudhe, N.} (2023)
Diffusion of binary opinions in a growing population with heterogeneous behaviour and external influence. Netw. Heterog. Media, \textbf{18}(3): 1288--1312.


\bibitem{NySz}
{\sc Nyczka, P. and Sznajd-Weron, K.} (2013)
Anticonformity or Independence?-Insights from Statistical Physics.
{\em J. Stat. Phys.} \textbf{151}(1-2), 174--202.

\bibitem{Peng2021}
{\sc Peng, K., Lu, Z., Lin, V., Lindstrom, M.R., Parkinson, C., Wang, C., Bertozzi, A.L., and Porter, M.A.} (2021)
A multilayer network model of the coevolution of the spread of a disease and competing opinions,
\textit{Math. Models Methods Appl. Sci.}, \textbf{31}(12), 2455--2494.

\bibitem{Qin} {\sc Qin, S.} (2025) Recurrence and transience of multidimensional elephant random walks. Ann. Probab. \textbf{53}(3): 1049--1078.

%\bibitem{stout} {\sc Stout, W.F.} (1974) {\it Almost sure convergence, Probability and Mathematical Statistics.} Vol. 24, Academic Press, New York-London.

\bibitem{Szna} {\sc Sznajd-Weron, K. and Sznajd, J.} (2000)
Opinion evolution in closed community, Int. J. Mod. Phys. C \textbf{11}: 1157–1165


\bibitem{gue}
{\sc Vázquez Guevara, V. H.} (2019)
On the almost sure central limit theorem for the elephant random walk, J. Phys. A: Math. Theor. \textbf{52}(1), 475201.

%\bibitem{Ze}
%{\sc Zelterman, D.} (2004)
%{\em Discrete Distributions. Applications in the Health Sciences}. John Wiley \& Sons, New York.

%\bibitem{ZH}
%{\sc Zhu, H. and Huberman, B. A.} (2014) To switch or not to switch understanding social influence in online choices. {\em Am. Behav. Sci.} \textbf{58}(10), 1329--1344.


\end{thebibliography}
\end{document}